\documentclass[a4paper,11pt, reqno]{amsart}
\usepackage{amssymb,amsthm,amsmath}
\usepackage{amsmath}
 \usepackage{enumerate}
\usepackage{ifthen}
\usepackage{graphicx}
\usepackage{cite}
\usepackage{tikz} 
\usetikzlibrary{arrows.meta}
\usepackage{comment}
\usepackage{amsmath,amscd,amssymb}
\usepackage{latexsym}
\usepackage[colorlinks,citecolor=blue,pagebackref,hypertexnames=false]{hyperref}

\numberwithin{equation}{section}
\allowdisplaybreaks[2]
\theoremstyle{plain}
\newtheorem{theorem}{Theorem}[section]

\newtheorem{lemma}[theorem]{Lemma}

\newtheorem{proposition}[theorem]{Proposition}

\theoremstyle{definition}

\theoremstyle{remark}

\newtheorem{case[theorem]}{Case}

\def \R{{\mathbb R}}

\def\norm#1.#2.{\lVert#1\rVert_{#2}}

\def\R{\mathbb R}

\title[Strichartz Estimates for the Schr\"odinger equation associated  to  OU operator]{ Local Dispersive and Strichartz estimates for the Schr\"odinger equation associated  to the Ornstein-Uhlenbeck operator}

\author{Aparajita Dasgupta}
\author{Uttam Kumar Dolai}
\author{Cheng Luo} 
\author{Manli Song}
\thanks{Corresponding author: Manli Song}

\address{\endgraf Department of Mathematics
Indian Institute of Technology, Delhi, Hauz Khas
New Delhi-110016
India}
\email{adasgupta@maths.iitd.ac.in}
\address{\endgraf Department of Mathematics
Indian Institute of Technology, Delhi, Hauz Khas
New Delhi-110016
India}
\email{maz248176@maths.iitd.ac.in}

\address{\endgraf School of Mathematics and Statistics, Northwestern Polytechnical University, Xi'an, Shaanxi 710129, China}
\email{caashmale@163.com}

\address{\endgraf School of Mathematics and Statistics, Northwestern Polytechnical University, Xi'an, Shaanxi 710129, China}
\email{mlsong@nwpu.edu.cn}


\date{\today}
\begin{document}

\begin{abstract}
In this paper we study the linear and nonlinear Schrödinger equations associated with the Ornstein-Uhlenbeck (OU) operator endowed with the Gaussian measure. While classical Strichartz estimates are well-developed for the free Schrödinger operator on Euclidean spaces, extending them to non-translation-invariant operators like the OU operator presents significant challenges due to the lack of global dispersive decay. In this work, we overcome these difficulties by deriving localized $L^1 \to L^\infty$ dispersive estimates for the OU Schrödinger propagator using Mehler kernel techniques. We then establish a family of weighted Strichartz estimates in Gaussian $L^p$ spaces via interpolation and the abstract $TT^*$-method. As an application, we prove local well-posedness results for the nonlinear Schrödinger equation with power-type nonlinearity in both subcritical and critical regimes. Our framework reveals new dispersive phenomena in the context of the OU semigroup and provides the first comprehensive Strichartz theory in this setting.
\end{abstract}

	\maketitle

	\allowdisplaybreaks

	\tableofcontents

	\section{Introduction }

    The Schr\"{o}dinger equation is a type of dispersive equation, which describes wave propagation with dispersion effects. Commonly observed in light waves, sound waves, and water waves, dispersion effect, informally, refers to the phenomenon where different frequencies of a wave propagate at different velocities, causing the wave to spread out over time. Strichartz estimates are fundamental tools in the study of dispersive partial differential equations, particularly in establishing well-posedness and scattering for nonlinear Schrödinger equations (NLS). While such estimates are classical in the Euclidean setting, recent decades have seen growing interest in extending them to more general operators and geometries, including those with underlying symmetry or weighted measures.

    Let us consider the linear Schrödinger equation posed on the real line:
\begin{equation}
\begin{cases}
i u_{t} - \Delta u = 0, & \; t \in \mathbb{R}, \\
u(0,x) = u_0(x), & x \in \mathbb{R}.
\end{cases}
\end{equation}
The solution operator associated with this equation, denoted by the unitary group \( (e^{it\Delta})_{t \in \mathbb{R}} \), satisfies two fundamental properties that can be derived, for example, using Fourier analysis (see \cite[Theorem IX.30]{MR}).

\medskip
First, we have the conservation of the \( L^2 \)-norm:
\[
\| e^{-it\Delta} u_0 \|_{L^2(\mathbb{R})} = \| u_0 \|_{L^2(\mathbb{R})}, \quad \text{for all } t \in \mathbb{R}.
\]

Second, the group exhibits a dispersive decay estimate:
\[
\| e^{-it\Delta} u_0 \|_{L^\infty(\mathbb{R})} \lesssim \frac{1}{\sqrt{|t|}} \| u_0 \|_{L^1(\mathbb{R})}, \quad \text{for all } t \ne 0.
\]

\medskip
By interpolating between the above two estimates, one obtains a family of time-decay bounds for intermediate \( L^p \)-spaces, namely:

\begin{itemize}
\item[(i)] \textbf{Dispersive estimates:} For any \( p \in [1,2] \),
\[
\| e^{-it\Delta} u_0 \|_{L^{p'}(\mathbb{R})} \lesssim |t|^{-\left(\frac{1}{2} - \frac{1}{p} \right)} \| u_0 \|_{L^p(\mathbb{R})}, \quad t \ne 0.
\]
\end{itemize}

\medskip
In addition to these pointwise-in-time estimates, a broader framework of \textit{space-time} estimates—known as \textbf{Strichartz estimates}—has been developed. These were established in a general setting by Keel and Tao \cite{KT}, and include the following types:

\begin{itemize}
\item[(ii)] \textbf{Strichartz estimates (homogeneous, dual homogeneous, and inhomogeneous):}
\begin{align*}
\| e^{-it\Delta} u_0 \|_{L_t^q(\mathbb{R}, L_x^r(\mathbb{R}))} &\lesssim \| u_0 \|_{L^2(\mathbb{R})}, \\
\left\| \int_{\mathbb{R}} e^{is\Delta} F(s, \cdot) \, ds \right\|_{L^2(\mathbb{R})} &\lesssim \| F \|_{L_t^{q'}(\mathbb{R}, L_x^{r'}(\mathbb{R}))}, \\
\left\| \int_{s < t} e^{-i(t-s)\Delta} F(s, \cdot) \, ds \right\|_{L_t^q(\mathbb{R}, L_x^r(\mathbb{R}))} &\lesssim \| F \|_{L_t^{\tilde{q}'}(\mathbb{R}, L_x^{\tilde{r}'}(\mathbb{R}))},
\end{align*}
\end{itemize}
where the pairs \( (q, r) \) and \( (\tilde{q}, \tilde{r}) \) are sharp \( 1/2 \)-admissible exponents. Here for any $\sigma>0$, we say that a pair $(q,r)$ is sharp $\sigma$-admissible if $q,r\geq2$, $(q,r,\sigma)\neq (2,\infty,1)$, and the following scaling condition holds:
\begin{equation}
\frac{1}{q}+\frac{\sigma}{r}=\frac{\sigma}{2}.
\end{equation}

  Strichartz estimates have become indispensable tools in the analysis of dispersive partial differential equations, especially in the study of well-posedness and long-time behavior of nonlinear Schrödinger equations (NLS). Although these estimates are classical for the free Schrödinger operator on $\mathbb{R}^d$, significant efforts have been made to extend them to more general geometric and algebraic settings, such as nilpotent Lie groups of step two (see \cite{BGX2000, H2005, BKG, Song2016, SZ, BBG2021, SY2023,FMV1, FV, LS2014, Song-Tan}), metric measure spaces (see \cite{LMS, FS, Ratna3, BDDM2019, DPR2010, FSW, NR2005, R2008, MS, FMSW}), hyperbolic spaces (see \cite{AP2009, AP2014, APV2012}), compact Riemannian manifolds (see \cite{B1993, burq2004strichartz}) and bounded domains (see \cite{ILP2014}), etc, each adapted to the symmetry and measure structure of the underlying space.  We also  refer to the monographs \cite{W. Urbina-Romero} and the references therein for an overview of this topic in the Euclidean framework and other settings.

   In the Euclidean setting, well-posedness theory, establishing the existence, uniqueness, and continuous dependence of solutions on initial data, has been extensively developed for various nonlinearities and function spaces for the nonlinear Schrödinger equation (NLS) with power nonlinearity.

 The foundational framework was established through Strichartz estimates, first introduced by Strichartz~\cite{R.S.Strichartz} in 1977 and refined in works such as Ginibre and Velo~\cite{ginibre1992smoothing} and Cazenave~\cite{TC}. These works showed that solutions to NLS equations with power-type nonlinearity
\[
i \partial_t u - \Delta u = \mu |u|^{p-1} u,
\]
are locally well-posed in $H^s(\mathbb{R}^d)$ for subcritical exponents $p < 1 + \frac{4}{d-2s}$. The key analytical tool is the dispersive decay of the linear Schrödinger group, allowing for contraction mapping arguments in Strichartz-type Banach spaces.
 In the critical setting---when the scaling symmetry of the NLS preserves the norm of the initial data---well-posedness becomes more delicate. Techniques such as concentration-compactness and profile decompositions have proven instrumental. The work of Kenig and Merle~\cite{kenig2006global} is a landmark result in this area, proving global well-posedness and scattering in the energy-critical focusing case with radial data. Subsequent works by Visan, Dodson, and others have extended these results to non-radial cases and other critical regimes.
 On non-Euclidean domains such as compact Riemannian manifolds, the dispersive properties of the Schrödinger flow are affected by the geometry. In such contexts, Burq, Gérard, and Tzvetkov~\cite{burq2004strichartz} demonstrated that Strichartz estimates generally incur a loss of derivatives, which restricts well-posedness to higher regularity spaces. For tori and product manifolds, additional refinements such as bilinear estimates and Bourgain-type spaces have been introduced~\cite{ionescu2009semilinear}.
 In recent years, significant interest has developed in the study of NLS equations involving nonstandard operators such as the Hermite operator, Laguerre operator, and Dunkl Laplacian. These arise naturally in settings with weighted measures or radial symmetry. The analysis often requires spectral and semigroup methods adapted to the operator structure, as seen in works by Thangavelu~\cite{thangavelu1993lectures} and Banica~\cite{banica2000remarks}. 
 
  In this work, we establish new Strichartz estimates for the Schrödinger equation associated with the Ornstein-Uhlenbeck (OU) operator
$$\mathbf{L}:=-\frac{1}{2}\Delta+x\cdot\nabla=\frac{1}{2}\nabla^*\nabla,$$
	where $\Delta=\sum_{k=1}^{d}\frac{\partial^2}{\partial x_k^2}$ is the Laplacian and $\nabla^*$ is the adjoint operator of the gradient operator $\nabla=\left(\frac{\partial}{\partial x_1},\frac{\partial}{\partial x_2},\cdots, \frac{\partial}{\partial x_d}\right)$ endowed with the standard Gaussian measure, $\gamma_d(x)dx=\pi^{-d/2}e^{-|x|^2}dx$. The OU operator is of central importance in stochastic analysis, semigroup theory, and infinite-dimensional analysis, yet its dispersive properties in the context of PDEs are poorly understood.  The study of Ornstein-Uhlenbeck operator experienced a lot of developments in the last couple of decades. For more details, we refer the reader to the book of Urbina~\cite{W. Urbina-Romero}. Our study develops a weighted Strichartz framework tailored to the OU operator, leveraging its explicit spectral decomposition in terms of Hermite polynomials and the Mehler kernel. The key contributions of this work are as follows:

\begin{itemize}
    \item We derive a \textbf{localized} $L^1 \to L^\infty$ \textbf{dispersive estimate} for the OU Schrödinger propagator $e^{itL}$, valid for $0 < |t| \leq \pi/2$, using the kernel formula via Mehler's expansion.
    
    \item Based on this, we establish a family of {weighted Strichartz estimates} in Gaussian $L^p$-spaces $L^p_{\gamma_d}(w)$, where the weight $w(x) = e^{-|x|^2/2}$ reflects the geometry of the OU semigroup.
    
    \item We apply these estimates to prove \textbf{local well-posedness} of the nonlinear Schrödinger equation
    \[
    i\partial_t u - Lu = \mu w^{p}|u|^{p-1} u, \quad u(0,x) = u_0(x),
    \]
    in the Gaussian $L^2$ space, both in the \textbf{subcritical regime} $1 < p < 1 + \frac{4}{d}$ and the \textbf{critical case} $p = 1 + \frac{4}{d}$. The analysis makes essential use of a Strichartz-type function space that incorporates both time integrability and spatial Gaussian weights.
\end{itemize}

Our methods combine harmonic analysis on Gaussian spaces, spectral theory, and interpolation techniques. In particular, the {abstract $TT^*$-method of Keel and Tao} \cite{KT} is adapted to the weighted setting, and the interpolation framework of Bergh and L\"{o}fstr\"om \cite{Bergh-Lofstrom} is used to characterize the involved function spaces. While Hermite and Laguerre operators exhibit oscillatory semigroups that admit dispersive estimates under careful analysis, the OU operator presents new challenges due to its lack of translation symmetry, non-unitarity of the flow in $L^p$, and kernel periodicity. These issues prevent global dispersive estimates and complicate any direct extension of the Euclidean theory. By focusing on local-in-time behavior and incorporating weight structures adapted to the Gaussian measure, we overcome these obstacles and provide the first rigorous dispersive framework for the OU Schr\"{o}dinger equation.

The primary novelty of this work lies in the derivation of weighted Strichartz estimates for the OU Schrödinger group $e^{it\mathbf{L}}$, which, to the best of our knowledge, has not been comprehensively developed in prior literature. Unlike the translation-invariant setting of the classical Laplacian, the OU semigroup lacks
global dispersive decay. We overcome this by deriving localized $L^{1}-L^{\infty}$ dispersive estimates in
Gaussian-weighted spaces and using interpolation theory (notably the Keel-Tao framework) to construct a family of admissible Strichartz estimates.

 Our second major contribution is to apply these weighted Strichartz estimates to prove local well-posedness of the nonlinear Schr\"{o}dinger equation with power-type nonlinearities, both in the subcritical
and critical regimes. We demonstrate that, despite the operator's lack of translation invariance and
compact resolvent, one can still establish robust control over the nonlinear evolution via a careful
functional-analytic framework in weighted spaces.

The structure of the paper is as follows. In Section 2, we review the spectral properties of the OU operator and its associated semigroup. Section 3 develops the dispersive estimates for both the classical and OU Schrödinger equations and weighted Strichartz estimates via interpolation and duality arguments. Section 4 applies these estimates to prove local well-posedness results for the nonlinear Schrödinger equation in both subcritical and critical regimes.

    \section{Preliminary spectral theory for the Ornstein-Uhlenbeck operator }
    In this section, we specify the definition of the Ornstein-Uhlenbeck operator $\mathbf{L}$ and recall some known facts about its spectrum. \\
    
	For any $k\in\mathbb{N}$, let $H_k$ denote the one-dimensional Hermite polynomial on $\mathbb{R}$ defined by
	$$H_k(x)=(-1)^k\frac{d^k}{dx^k}(e^{-x^2})e^{x^2},\forall x\in \mathbb{R},$$
	and $h_k$ be the one-dimensional normalized Hermite polynomial on $\mathbb{R}$ defined by
	$$h_k(x)=(2^kk!)^{-\frac{1}{2}}H_k(x).$$
	The Ornstein-Uhlenbeck operator with the Gaussian measure on $\mathbb{R}$ is given by
	$$L:=-\frac{1}{2}\frac{\partial^2}{\partial x^2}+x\frac{\partial}{\partial x}.$$
	Then the one-dimensional normalized Hermite polynomials $\{h_k\}_{k=0}^\infty$ which form a complete orthonormal system in the Gaussian $L^2$ space $L^2(\mathbb{R},\gamma_1(x)dx)$, where $\gamma_1(x)=\pi^{-\frac{1}{2}}e^{-|x|^2}dx$ is the Gaussian measure on $\mathbb{R}$. Eigenfunctions of the Ornstein-Uhlenbeck operator $\mathbf{L}$ with the eigenvalues $k$ are the normalized Hermite polynomials, i.e.,
	$$Lh_k=kh_k.$$
	Therefore, for any $f\in L^2(\mathbb{R},\gamma_1(x)dx)$, it has the Hermite expansion
	$$f=\sum_{k=0}^\infty \langle f,h_k\rangle_{\gamma_1}h_k,$$
	where $\langle\cdot,\cdot\rangle_{\gamma_1}$ is an inner product inherited from $L^2(\mathbb{R},\gamma_1(x)dx)$.

	For any multi-index $\alpha=(\alpha_1,\alpha_2,\cdots,\alpha_d)\in\mathbb{N}^d$, $d\in\mathbb{N}^+$, the $d$-dimensional normalized Hermite polynomial $\mathbf{h}_\alpha$ is defined by the tensor product of one-dimensional normalized Hermite polynomials given by
	$$\mathbf{h}_\alpha(x)=\prod_{k=1}^{d}h_{\alpha_k}(x_k), \forall x=(x_1,x_2,\cdots,x_d)\in\mathbb{R}^d.$$
	Moreover, the $d$-dimensional Ornstein-Uhlenbeck operator associated with the Gaussian measure $$\gamma_d(x)dx:=\pi^{-\frac{d}{2}}e^{-|x|^2}dx_1dx_2\cdots dx_d$$ on the Euclidean space $\mathbb{R}^d$  is defined by
	$$\mathbf{L}:=\sum_{k=1}^{d}\left[-\frac{1}{2}\frac{\partial^2}{\partial x_k^2}+x_k\frac{\partial}{\partial x_k}\right]==-\frac{1}{2}\Delta+x\cdot\nabla=\frac{1}{2}\nabla^*\nabla,$$
	where $\Delta=\sum_{k=1}^{d}\frac{\partial^2}{\partial x_k^2}$ is the Laplacian and $\nabla^*$ is the adjoint operator of the gradient operator $\nabla=\left(\frac{\partial}{\partial x_1},\frac{\partial}{\partial x_2},\cdots, \frac{\partial}{\partial x_d}\right)$ with respect to the Gaussian measure $\gamma_d(x)dx$.
The family of $\{\mathbf{h}_\alpha\}_{\alpha\in\mathbb{N}^d}$ also forms an orthonormal system in the Gaussian $L^2$ space $L^2(\mathbb{R}^d,\gamma_d(x)dx)$. Also,  $\mathbf{h}_\alpha$s are the eigenfunctions of $\mathbf{L}$
with the eigenvalues $|\alpha|:=\sum\limits_{k=1}^d\alpha_k$, i.e.,
$$\mathbf{L}\mathbf{h}_\alpha=|\alpha|\mathbf{h}_\alpha.$$

For any given $f\in L^2(\mathbb{R}^d,\gamma_d(x)dx)$, it has the Hermite expansion
$$f=\sum_{\alpha\in\mathbb{N}^d}\langle f,\mathbf{h}_\alpha\rangle_{\gamma_d}\mathbf{h}_\alpha=\sum_{k=0}^\infty\mathbf{J}_kf,$$
where $\mathbf{J}_k$ denotes the orthogonal projection operator corresponding to the eigenvalue $k$, i.e.,
$$\mathbf{J}_kf=\sum\limits_{|\alpha|=k}\langle f,\mathbf{h}_\alpha\rangle_{\gamma_d}\mathbf{h}_\alpha.$$
	
The Ornstein-Uhlenbeck operator $\mathbf{L}$ generates a semigroup called the Ornstein-Uhlenbeck semigroup $e^{-t\mathbf{L}}, t>0$, defined by
$$e^{-t\mathbf{L}} f=\sum_{k=0}^\infty e^{-kt}\mathbf{J}_k f,$$
for $f\in L^2(\mathbb{R}^d,\gamma_d(x)dx)$. The above expression can also be written as
$$e^{-t\mathbf{L}} f(x)=\int_{\mathbb{R}^d} M_t(x,y)f(y)\gamma_d(y)dy,$$
where the kernel $M_t(x,y)$ is given by the expansion
$$M_t(x,y)=\sum_{\alpha\in\mathbb{N}^d} e^{-|\alpha|t}\mathbf{h}_\alpha(x)\mathbf{h}_\alpha(y).$$
By Mehler's formula $z=r+it, r>0,t\in\mathbb{R}$, the kernel of the operator $e^{-z\mathbf{L}},$  given by 
$$M_z(x,y)=\sum_{\alpha\in\mathbb{N}^d} e^{-|\alpha|z}\mathbf{h}_\alpha(x)\mathbf{h}_\alpha(y)=\sum_{k=0}^\infty \sum_{|\alpha|=k} e^{-kz}\mathbf{h}_\alpha(x)\mathbf{h}_\alpha(y),$$ can be written as
$$M_z(x,y)=\frac{e^{\frac{dz}{2}}}{(2 \sinh z)^\frac{n}{2}}e^{\frac{1}{2}\left((1-\coth z)(|x|^2+|y|^2)+\frac{2x\cdot y}{\sinh z}\right)}.$$
For $t\in \mathbb{R}\setminus \pi\mathbb{Z}$, letting $r\rightarrow 0^+$, the kernel of the operator $e^{-it\mathbf{L}}$ can also be written as
\begin{equation}\label{schrodinger-kernel}
M_{it}(x,y)=\frac{e^{-\frac{id\pi}{4}}e^{\frac{idt}{2}}}{(2 \sin t)^\frac{d}{2}}e^{\frac{1}{2}(|x|^2+|y|^2)}e^{\frac{i}{2}\left(\cot t(|x|^2+|y|^2)-\frac{2x\cdot y}{\sin t}\right)}.
\end{equation}
Moreover, we have
\begin{equation}\label{period}
M_{it}(x,y)=\overline{M_{-it}(x,y)} \text{ and } M_{i(t+\pi)}(x,y)=M_{it}(-x,y).
\end{equation}

In the sequel, we shall denote the Gaussian $L^p$ space on $\mathbb{R}^d$ by $$L_{\gamma_d}^p(\mathbb{R}^{d}):=L^p(\mathbb{R}^d,\gamma_d(x)dx),$$  and its weighted spaces by $$L_{\gamma_d}^p(w):=L^p(\mathbb{R}^d,w(x)\gamma_d(x)dx),$$ where $w$ is a positive even weight function on $\mathbb{R}^d$. By \eqref{period},the $L_{\gamma_d}^p(w)$ norm of $e^{-it\mathbf{L}}f$ is $\pi$-periodic in $t$ and thus determined by its values for $t\in [-\frac{\pi}{2},\frac{\pi}{2}]$.

Here we recall the definition of  the \textit{mixed Lebesgue space} \( L^q_t L^r_x(I \times \mathbb{R}^d) \) over a time interval \( I \subset \mathbb{R} \) as the space of all measurable functions \( f: I \times \mathbb{R}^d \to \mathbb{C} \) such that  
\[
\|f\|_{L^q_t L^r_x} := \left( \int_I \left( \int_{\mathbb{R}^d} |f(t,x)|^r\, \mathrm{d}x \right)^{\frac{q}{r}} \mathrm{d}t \right)^{\frac{1}{q}} < \infty,
\]
for \( 1 \leq q, r < \infty \). When \( q = \infty \) or \( r = \infty \), the norm is modified in the usual way. For example, if \( q = \infty \), we define
\begin{equation}
\label{eq:Linfty}
\|f\|_{L^{\infty}_t L^r_x} := \text{ess} \sup_{t \in I} \|f(t,\cdot)\|_{L^r_x(\mathbb{R}^d)}.
\end{equation}

We also consider the \textit{Gaussian-weighted mixed Lebesgue space} \( L^q_t(I; L^r_{\gamma_d}(w)) \), consisting of all measurable functions \( f: I \times \mathbb{R}^d \to \mathbb{C} \) for which
\[
\|f\|_{L^q_t(I; L^r_{\gamma_d}(w))} := \left( \int_I \left( \int_{\mathbb{R}^d} |f(t,x)|^r\, w(x)\, \gamma_d(x)\, \mathrm{d}x \right)^{\frac{q}{r}} \mathrm{d}t \right)^{\frac{1}{q}} < \infty,
\]
where \( \gamma_d(x) \) denotes the standard Gaussian density on \( \mathbb{R}^d \), and \( w(x) \) is a non-negative measurable weight function. In the case \( q = \infty \) or \( r = \infty \), the norm is defined analogously as in \eqref{eq:Linfty}.

\section{Strichartz estimates for the Ornstein Uhlenbeck operator}
In this section, we aim to derive one of the principal results of this work: the Strichartz estimates for the Ornstein-Uhlenbeck operator. As a preliminary step, we recall the classical Strichartz estimates associated with the Laplacian, which serve as a guiding framework for our analysis. The original formulation of Strichartz estimates in \( L^p \) spaces was introduced by R.~Strichartz~\cite{R.S.Strichartz}, in a context that was closely connected to the restriction problem for the Fourier transform on submanifolds. Subsequently, Ginibre and Velo~\cite{ginibre1992smoothing} employed the \( TT^* \) method to obtain such estimates, providing a powerful functional-analytic tool that has influenced much of the subsequent development in this area. 
As is well known, Fourier analysis plays a fundamental role in the study of the classical free Schr\"odinger equation on \( \mathbb{R}^d \) for \( d \geq 1 \). Consider 
\begin{align*}
		\begin{cases}
			i\partial_t u-\Delta u=0,\quad  t \in \mathbb{R}, \\
			u(0,x)=u_0(x)\in L^2(\mathbb{R}^d).
		\end{cases}
	\end{align*}
The solution can be written explicitly as
\begin{align}
\label{Homogeneoussolution}u(t,x)&=e^{-it\Delta}u_0(x)\\
      &\nonumber=\frac{1}{(2\pi)^d}\int_{\mathbb{R}^d} e^{it(|\xi|^2+x\cdot \xi)} \widehat{u_0}(\xi) \mathrm{d}\xi\\
      &\nonumber=\frac{1}{(4\pi i t)^\frac{d}{2}}\int_{\mathbb{R}^d} e^{i\frac{|x-y|^2}{4t}} u_0(y) \mathrm{d} y.
\end{align}
 Therefore, we have the following dispersive estimate $L^1\to L^\infty$.
\begin{equation}\label{disest}
\|u(t,\cdot)\|_{L^\infty(\mathbb{R}^d)}=\|e^{it\Delta}u_0\|_{L^\infty(\mathbb{R}^d)}\lesssim \frac{1}{|t|^\frac{d}{2}}\|u_0\|_{L^1(\mathbb{R}^d)},\quad \forall t\ne 0.
\end{equation}
The corresponding inhomogeneous equation is 
\begin{align}
\label{Laplacian NLS}
\begin{cases}
i\partial_t u-\Delta u=F(t,x),\quad x \in \mathbb{R}^d, t \in \mathbb{R}, \\
u(0,x)=u_{0}(x).
\end{cases}
\end{align}
By Duhamel's principle, the integral version of \eqref{Laplacian NLS} has the form 
\begin{align}
\label{Laplacian Duhamel}
u(t,x)=e^{-it\Delta}u_0(x)-i\int_{0}^{t}e^{-i(t-s)\Delta}F(s,\cdot)\ \mathrm{d}s.
\end{align}

From the above dispersive estimate (\ref{disest}), using the $TT^*$ argument, the Strichartz estimate for the free Schr\"odinger equation can be deduced. Here we recall the abstract Stricharz estimates established by Keel-Tao \cite{KT}.
\begin{lemma} Let $(X,dx)$ be a measure space and $H$ be a Hilbert space. Suppose that for each time $t\in\mathbb{R}$ we have an operator $U(t): H\to L^2(X)$ which obey the energy estimate, i.e., for all  we have
\begin{equation*}
\|U(t)\|_{L^2(X)}\lesssim \|f\|_H,\;\forall t\in\mathbb{R},\;f\in H,
\end{equation*}
and that for some $\sigma>0$ the following decay estimate holds
\begin{equation*}
\|U(s)(U(t))^*g\|_{L^\infty(X)}\lesssim |t-s|^{-\sigma} \|g\|_{L^1(X)},\;\forall t\neq s\in\mathbb{R},\; g\in L^1(X).
\end{equation*}
Then we have the estimates 
\begin{align*}
\|U(t)f\|_{L^{q}_{t}(\mathbb{R},L^r(X))}&\lesssim \|f\|_H,\\
\left\|\int_{\mathbb{R}}U(t)(U(s))^*F(s,\cdot)\right\|_{L_{t}^{q}(\mathbb{R},L^r(X))}&\lesssim \|F\|_{L_{t}^{\tilde{q}'}(\mathbb{R},L^{\tilde{r}'}(X))},
\end{align*}
hold for all sharp $\sigma$-admissible pairs $(q,r),(\tilde{q},\tilde{r})$.
\end{lemma}
As a consequence of the above theorem, taking $U(t)=e^{-it\Delta}$ and $H=L^2(\mathbb{R}^d)$, one can obtain the Strichartz estimates for the classical Schr\"odinger operator.
\begin{theorem}
Suppose that $d\geq 1$ and  $(q,r),(\tilde{q},\tilde{r})$ are sharp $\frac{d}{2}$-admissible pairs. Then
\begin{align*}
\|e^{-it\Delta}u_0\|_{L_{t}^{q}(\mathbb{R},L^r(\mathbb{R}^d))}&\lesssim \|u_0\|_{L^2(\mathbb{R}^d)},\\
\left\|\int_{\mathbb{R}}e^{-i(t-s)\Delta}F(s,\cdot)\right\|_{L_{t}^{q}(\mathbb{R},L^r(\mathbb{R}^d))}&\lesssim \|F\|_{L_{t}^{\tilde{q}'}(\mathbb{R},L^{\tilde{r}'}(\mathbb{R}^d))}.
\end{align*}
\end{theorem}

Now let us consider the free Schr\"odinger equation related to the Ornstein-Uhlenbeck Operator on $\mathbb{R}^d$ 
\begin{align}\label{eq-OU}
		\left\{\begin{array}{l}
			i\partial_t u- \mathbf{L} u=0,\quad  t \in \mathbb{R}, \\
			u(0,x)=u_0(x) \quad x \in \mathbb{R}^d.
		\end{array}\right.
	\end{align}
The solution of \eqref{eq-OU} can be written explicitly by
\begin{align*}
u(t,x)&=e^{-it\mathbf{L}}u_0(x)\\
      &=\int_{\mathbb{R}^d} M_{it}(x,y) u_0(y) \gamma_d(y)dy.
\end{align*}
It follows from \eqref{schrodinger-kernel} that
\begin{equation*}
|u(t,x)|=|e^{-it\mathbf{L}}u_0(x)|\lesssim \frac{1}{|\sin t|^\frac{d}{2}}\int_{\mathbb{R}^d}e^{\frac{1}{2}(|x|^2+|y|^2)}|u_0(y)|\gamma_d(y)dy,\quad \forall t\in \mathbb{R}\setminus \pi\mathbb{Z}.
\end{equation*}
Let the weight function $w(x)=e^{-\frac{|x|^2}{2}}$ and we have
\begin{equation*}
|u(t,x)w(x)|=|e^{-it\mathbf{L}}u_0(x)w(x)|\lesssim \frac{1}{|\sin t|^\frac{d}{2}}\int_{\mathbb{R}^d}|u_0(y)|w^{-1}(y)\gamma_d(y)dy,\quad \forall t\in \mathbb{R}\setminus \pi\mathbb{Z}.
\end{equation*}
Hence, for all $ 0<|t|\leq \frac{\pi}{2}$ we obtain the local weighted $L^1\to L^\infty$ dispersive estimate as
\begin{equation}\label{weighted dispersive}
\|u(t,\cdot)\|_{L_{\gamma_d}^\infty(w)}=\|e^{-it\mathbf{L}}u_0\|_{L_{\gamma_d}^\infty(w)}\lesssim \frac{1}{|\sin t|^\frac{d}{2}}\|u_0\|_{L_{\gamma_d}^1(w^{-1})}\lesssim \frac{1}{| t|^\frac{d}{2}}\|u_0\|_{L_{\gamma_d}^1(w^{-1})}.
\end{equation}

Using the above dispersive estimate \eqref{weighted dispersive}, we can derive the following weighted Strichartz estimates for the solution of \eqref{eq-OU}.
\begin{theorem}\label{main}
Let $d\geq1$ and $(q,r)$, $(\tilde{q},\tilde{r})$ be sharp $\frac{d}{2}$-admissible pairs. If $d=1$ and $(q,r)=(4,\infty)$, the Strichartz estimates hold
\begin{equation*}
\|e^{-it\mathbf{L}}u_0\|_{L^{4}_{t}([-\frac{\pi}{2},\frac{\pi}{2}],L^\infty_{\gamma_1}(w))}\lesssim \|u_0\|_{L^2_{\gamma_1}}.
\end{equation*}
For other sharp $\frac{d}{2}$-admissible pairs $(q,r)$, one has the weighted Strichartz estimates 
\begin{equation}
\label{HomogeneousStrichartz}
\|e^{-it\mathbf{L}}u_0\|_{L^{q}_{t}([-\frac{\pi}{2},\frac{\pi}{2}],L^r_{\gamma_d}(w^{r-2}))}\lesssim \|u_0\|_{L^2_{\gamma_d}},
\end{equation}
and
\begin{equation}
\label{HomogeneousDualStrichartz}
\bigg\|\int_{-\frac{\pi}{2}}^{\frac{\pi}{2}}e^{i(t-s)\mathbf{L}}F(s,\cdot)\ \mathrm{d}s\bigg\|_{L_{t}^{q}([-\frac{\pi}{2},\frac{\pi}{2}],L^{r}_{\gamma_d}(w^{r-2}))}\lesssim \|F\|_{L_{t}^{\tilde{q}'}([-\frac{\pi}{2},\frac{\pi}{2}], L^{\tilde{r}'}_{\gamma_{d}}(w^{\tilde{r}'-2}))},
\end{equation}
where the weight function is 
\begin{equation*}
w(x)=e^{-\frac{|x|^2}{2}}.
\end{equation*}
\end{theorem}
In order to prove the above theorem, we will apply the general lemma of Keel-Tao \cite{KT}.
\begin{lemma}\label{abstract KT}
Let $\sigma>0$, $H$ be a Hilbert space and $B_0,B_1$ be Banach spaces. Suppose that for each time $t\in I$ we have an operator $U(t): H\to B_0^*$ such that
\begin{align*}
\|U(t)\|_{H\to B_0^*}&\lesssim 1,\\
\|U(t)(U(s))^*\|_{B_1\to B_1^*}&\lesssim |t-s|^{-\sigma}.
\end{align*}
Let $B_\theta$ denote the real interpolation space $(B_0,B_1)_{\theta,2}$. Then we have the estimates
\begin{align}\label{non}
\|U(t)f\|_{L^{q}_{t}(I,B^*_\theta)}&\lesssim \|f\|_H,\\
\label{otherdualpart}\left\|\int_I U(t)(U(s))^*F(s,\cdot)ds\right\|_{L^{q}_{t}(I,B^*_\theta)}&\lesssim \|F\|_{L_{t}^{\tilde{q}'}(I,B_{\tilde{\theta}})},
\end{align}
where $0\leq\theta\leq 1$, $2\leq q=\frac{2}{\sigma\theta}$, $(q,\theta,\sigma)\neq (2,1,1)$ and similarly for $(\tilde{q},\tilde{\theta})$. 
\end{lemma}

\begin{proof}[Proof of Theorem \ref{main}] In our case, let $U(t)=e^{-it\mathbf{L}}$ and
\begin{align*}
H=B_0&=L^2_{\gamma_d},\;B_1=L^1_{\gamma_d}(w^{-1}),\;\sigma=\frac{d}{2},\;2\leq \frac{4}{d\theta},\\
&0\leq \theta\leq 1,\; \left(\frac{4}{d\theta},\theta,\sigma\right)\neq (2,1,1).
\end{align*}
Applying the lemma \ref{abstract KT}, we get
\begin{equation*}
\|e^{-it\mathbf{L}}u_0\|_{L_{t}^\frac{4}{d\theta}([-\frac{\pi}{2},\frac{\pi}{2}],(L^2_{\gamma_d}, L^1_{\gamma_d}(w^{-1}))^*_{\theta,2})}\lesssim \|u_0\|_{L^2_{\gamma_d}}.
\end{equation*}
By Theorem 3.4.1. of Bergh-L\"ofstr\"om \cite{Bergh-Lofstrom}, for any $1\leq h\leq 2$, 
\begin{equation*}
(L^2_{\gamma_d}, L^1_{\gamma_d}(w^{-1}))_{\theta,h}\subseteq(L^2_{\gamma_d}, L^1_{\gamma_d}(w^{-1}))_{\theta,2},
\end{equation*}
and hence
\begin{equation}\label{interpolation}
\|e^{-it\mathbf{L}}u_0\|_{L_{t}^{\frac{4}{d\theta}}([-\frac{\pi}{2},\frac{\pi}{2}],(L^2_{\gamma_d}, L^1_{\gamma_d}(w^{-1}))^*_{\theta,h})}\lesssim \|u_0\|_{L^2_{\gamma_d}}.
\end{equation}
The conditions on $\theta$ can be divided into two cases:
\begin{enumerate}
\item when $d\geq 2$, $0\leq \theta\leq \frac{2}{d}$ and $(\theta,d)\neq (1,2)$;
\item when $d=1$, $0\leq \theta\leq 1$.
\end{enumerate}

\textbf{Case 1: } Let $d\geq 2$, $0\leq \theta\leq \frac{2}{d}$ and $(\theta,d)\neq (1,2)$. In this case, $0\leq \theta<1$. Indeed, for the value $\theta=0$, the operator $e^{it\mathbf{L}}$ obeys the energy conservation, which implies the Strichartz estimates.
\begin{equation}\label{energy conservation}
\|e^{-it\mathbf{L}}u_0\|_{L_{t}^{\infty}([-\frac{\pi}{2},\frac{\pi}{2}],L^2_{\gamma_d})}\lesssim \|u_0\|_{L^2_{\gamma_d}}.
\end{equation}
For any $0<\theta<1$, by Theorem 5.5.1. of Bergh-L\"ofstrom \cite{Bergh-Lofstrom}, 
\begin{equation*}
(L^2_{\gamma_d}, L^1_{\gamma_d}(w^{-1}))_{\theta,h}=L^h_{\gamma_d}(\tilde{w}),
\end{equation*}
with
\begin{equation*}
\frac{1}{h}=\frac{1-\theta}{2}+\frac{\theta}{1} \text{ and } \tilde{w}=w^{-h\theta},
\end{equation*}
which indicates that $1<h=\frac{2}{1+\theta}<2$. Therefore,
\begin{equation*}
(L^2_{\gamma_d}, L^1_{\gamma_d}(w^{-1}))_{\theta,h}=L^\frac{2}{1+\theta}_{\gamma_d}(w^{-\frac{2\theta}{1+\theta}}).
\end{equation*}
From the classical formula
\begin{equation*}
(L^\alpha_{\gamma_d}(w))^*=L^{\alpha'}_{\gamma_d}(w^{-\frac{\alpha'}{\alpha}}),\;\forall 1<\alpha<\infty.
\end{equation*}
It follows that
\begin{equation*}
(L^2_{\gamma_d}, L^1_{\gamma_d}(w^{-1}))_{\theta,h}^*=L^\frac{2}{1-\theta}_{\gamma_d}(w^{\frac{2\theta}{1-\theta}}).
\end{equation*}

Therefore, for $0<\theta<1$, the inequality \eqref{interpolation} gives the weighted Strichartz estimates.
\begin{equation}\label{weighted Strichartz}
\|e^{-it\mathbf{L}}u_0\|_{L^\frac{4}{d\theta}_{t}([-\frac{\pi}{2},\frac{\pi}{2}],L^\frac{2}{1-\theta}_{\gamma_d}(w^{\frac{2\theta}{1-\theta}}))}\lesssim \|u_0\|_{L^2_{\gamma_d}},
\end{equation}
which also holds for $\theta=0$ from \eqref{energy conservation}. 

\textbf{Case 2: } Let $d=1$, $0\leq \theta\leq 1$. For $0\leq \theta<1$. Proceeding in a manner similar to Case 1, one can derive the weighted Strichartz estimates as well, \eqref{weighted Strichartz}. For $\theta=1$, from \eqref{non} it follows that
\begin{equation*}
\|e^{-it\mathbf{L}}u_0\|_{L_{t}^{4}([-\frac{\pi}{2},\frac{\pi}{2}],L^\infty_{\gamma_1}(w))}=\|e^{-it\mathbf{L}}u_0\|_{L_{t}^{4}\left([-\frac{\pi}{2},\frac{\pi}{2}],\left(L^1_{\gamma_1}(w^{-1})\right)^*\right)}\lesssim \|u_0\|_{L^2_{\gamma_1}}.
\end{equation*}

For $0\leq \theta<1$, denoting $q=\frac{4}{d\theta}$ and $r=\frac{2}{1-\theta}$, we have $\frac{2\theta}{1-\theta}=-2+\frac{2}{1-\theta}=r-2$. The pair $(q,r)$ satisfies $\frac{2}{q}+\frac{\sigma}{r}=\frac{\sigma}{2}$ with $q\geq2$, $2\leq r<\infty$, $\sigma=\frac{d}{2}$ and $(q,r,d)\neq (2,\infty,2)$.
Therefore, we have the weighted Strichartz estimates
\begin{equation}\label{weighted Strichartz2}
\|e^{-it\mathbf{L}}u_0\|_{L_{t}^{q}([-\frac{\pi}{2},\frac{\pi}{2}],L^r_{\gamma_d}(w^{r-2}))}\lesssim \|u_0\|_{L^2_{\gamma_d}}.
\end{equation}

Now, using duality arguments  and proceeding with analogous computations as above, one can establish the homogeneous dual Strichartz estimate \eqref{HomogeneousDualStrichartz}.

\end{proof}

    \section{The well-posedness results}

In this section, we establish well-posedness results for a class of nonlinear Schr\"odinger equations (NLS) as an application of the Strichartz estimates derived earlier. We consider the Cauchy problem associated with the Ornstein–Uhlenbeck operator:
\begin{align}
\label{NLS}
\begin{cases}
i\partial_{t}u - \mathbf{L}u = F(u),\quad  t \in \mathbb{R}, \\
u(x,0) = u_{0}(x),\quad  x\in\mathbb{R}^{d},
\end{cases}
\end{align}
where the nonlinearity is of power type, given by \( F(u) = \mu w^p |u|^{p-1}u \), with \( p > 1 \) and \( \mu = \pm 1 \). The sign of \( \mu \) distinguishes the defocusing case (\( \mu = +1 \)) from the focusing case (\( \mu = -1 \)). We assume that the initial data \( u_0 \) belongs to the Gaussian space \( L^{2}_{\gamma_{d}} (\mathbb{R}^{d})\). In this setting, the critical exponent corresponds to \( p = 1 + \frac{4}{d} \). We will establish local well-posedness results for both the subcritical and critical cases.

The classical theory for nonlinear Schrödinger equations (NLS) typically applies when the initial data is very smooth and the nonlinearity \( u \mapsto F(u) \) is sufficiently regular. However, for rougher initial data or less regular nonlinearities, one must adopt a different framework. Since the differential formulation of the NLS,
\begin{equation*}
i\partial_t u - \mathbf{L} u = F(u),
\end{equation*}
requires a high degree of differentiability, it is often more effective to work instead with its integral (or Duhamel) formulation:
\begin{equation}
\label{Duhamel}
u(t,x) = e^{-it\mathbf{L}} u_0(x) - i \int_0^t e^{-i(t-s)\mathbf{L}} F(u(s,x))\, \mathrm{d}s.
\end{equation}

We may express this representation more concisely as
\[
u = u_0^* + D N(u),
\]
where \( u_0^*(t,x) := e^{-it\mathbf{L}} u_0(x) \) denotes the solution to the free Schrödinger equation,
\[
i\partial_t u + \mathbf{L} u = 0,
\]
\( N(u) = \mu w^p |u|^{p-1} u \) represents the nonlinearity, and \( D \) denotes the Duhamel operator defined by
\[
Dg(t,x) := -i\int_0^t e^{-i(t-s)\mathbf{L}} g(s,x)\, \mathrm{d}s.
\]

The integral equation \eqref{Duhamel} remains meaningful even when \( u \) is only a tempered distribution, provided that it belongs locally to the space \( L^q_t L^q_{\gamma_d}(I \times \mathbb{R}^d) \) for some \( q \geq 1 \). We refer to such functions \( u \) as distributional solutions to \eqref{NLS}. These solutions agree with classical ones when \( u \) is smooth, but they provide a broader framework that accommodates rough initial data.

In our setting, we consider initial data \( u_0 \in L^2_{\gamma_d}(\mathbb{R}^d) \). We say that a function \( u \) is a strong \( L^2_{\gamma_d}(\mathbb{R}^d) \) solution to \eqref{Duhamel} on a time interval \( I \) if:
\begin{itemize}
    \item \( u \) is a distributional solution to \eqref{Duhamel}, and
    \item \( u \in C^0_{t,\text{loc}}(I; L^2_{\gamma_d}(\mathbb{R}^d)) \).
\end{itemize}

A fundamental tool in well-posedness theory is the \textit{contraction theorem}. 
To establish the required contraction mapping, we rely on an abstract result formulated as Proposition 1.38 in \cite{T Tao}. The proof of this proposition can be found in \cite{E. Cordero3}.

\begin{proposition}[Iteration Principle]
\label{Iterationprinciple}
Let $\mathcal{ N }$ and \( \mathcal{S} \) be two Banach spaces. Let  $D : \mathcal{N} \to \mathcal{S}$ be a bounded  linear operator  with the bound 
\begin{align}
\label{Dandf}
\|Df\|_\mathcal{S} \leq C \|f\|_\mathcal{N}
\end{align}
for all \( f \in \mathcal{N} \) and some constant $C>0$ and let $N:\mathcal{S}\to\mathcal{N}$ with $N(0)=0$, be a nonlinear operator which is Lipschitz continuous and obeys the bounds
\begin{align}
\label{inequalityofN}
\|N(g) - N(g')\|_\mathcal{N} \leq \frac{1}{2C} \|g - g'\|_\mathcal{S}
\end{align}
for all \( g, g' \) in \( B_{\epsilon} = \{ g \in \mathcal{S} : \|g\|_\mathcal{S} \leq \epsilon \} \), for some \( \epsilon > 0 \). Then for all \( u_{0}^{*} \in B_{\frac{\epsilon}{2}} \),
there exists a unique solution \( u \in B_{\epsilon} \) to the equation
\begin{align}
\label{Fixtequatiion}
u = u_{0}^{*} + DN(u)
\end{align}
with Lipschitz map $u_{0}^{*}\mapsto u$ with constant at most $2$ i.e., if \( v \in B_{\epsilon} \) solves \( v = v_{0}^{*} + DN(v) \) for \( v_{0}^{*}\in B_{\frac{\epsilon}{2}} \), then
\begin{align}
\label{uvu*v*}
\|u - v\|_{\mathcal{S}(I\times\mathbb{R}^{d})} \leq 2 \|u_{0}^{*} - v_{0}^{*}\|_{\mathcal{S}}.
\end{align}
\end{proposition}

Proposition \ref{Iterationprinciple} serves as a fundamental tool in the analysis carried out in \cite{BArpad}, \cite{D'Ancona}, \cite{E. Cordero}, \cite{E. Cordero2}, and \cite{B. Wang}.

To investigate the well-posedness of the Cauchy problem, it is convenient to work within a function space \(\mathcal{S}(I \times \mathbb{R}^d)\) that simultaneously encodes all relevant Strichartz norms at a fixed level of regularity. In our setting, the initial data belong to the Gaussian-weighted space \(L^2_{\gamma_d}(\mathbb{R}^d)\). Accordingly, we define the following Strichartz space,

 \begin{align*}
 \mathcal{S}(I\times\mathbb{R}^{d})=\bigg\lbrace f: I\times\mathbb{R}^{d}\to\mathbb{C} : \sup\limits_{(q,r)\in\mathcal{A}}\|f\|_{L^{q}_{t}L^{r}_{\gamma_{d}}(w^{r-2})}<\infty \bigg\rbrace,
 \end{align*} 
which is a Banach space equipped with the norm 
\begin{align*}
\|f\|_{\mathcal{S}(I\times\mathbb{R}^{d})}=\sup\limits_{(q,r)\in\mathcal{A}}\|f\|_{L^{q}_{t}L^{r}_{\gamma_{d}}(w^{r-2})},
\end{align*}
where $\mathcal{A}$ is the set of all sharp $\frac{d}{2}$-admissible pairs. 
Throughout this section, we use the notation
\[
L^{q}_{t} L^{r}_{\gamma_{d}}(w^{r-2}) := L^{q}(I; L^{r}_{\gamma_{d}}(w^{r-2})),
\]
where the weight function is defined by \( w(x) = e^{-\frac{|x|^{2}}{2}} \), and the time interval \( I \subseteq \left[-\frac{\pi}{2}, \frac{\pi}{2}\right] \).

%
%
%
We define the space \(\mathcal{N}(I \times \mathbb{R}^d)\) as the dual space of \(\mathcal{S}(I \times \mathbb{R}^d)\).
If \( f \in \mathcal{N}(I \times \mathbb{R}^d) \), then for every admissible pair \((q,r) \in \mathcal{A}\), 
\begin{align*}
\|f\|_{\mathcal{N}}&=\sup_{\substack{\|g\|_{\mathcal{S}(I\times\mathbb{R}^{d})}\leq 1}}\left|\int_{I\times \mathbb{R}^{d}}f(t,x)g(t,x)\ \gamma_{d}(x)\mathrm{d}x \mathrm{d}t\right|\\
&\leq\sup_{\substack{\|g\|_{\mathcal{S}(I\times\mathbb{R}^{d})}\leq 1}}\left(\int_{I}\left(\int_{\mathbb{R}^{d}}|g(t,x)|^{r}w^{r-2}(x)\ \gamma_{d}(x)\mathrm{d}x\right)^{\frac{q}{r}}\ \mathrm{d}t\right)^{\frac{1}{q}}\\ 
&\times \left(\int_{I}\left(\int_{\mathbb{R}^{d}}|f(t,x)|^{r'}w^{-\frac{(r-2)r'}{r}}(x)\gamma_{d}(x)\ \mathrm{d}x\right)^{\frac{q'}{r'}}\ \mathrm{d}t\right)^{\frac{1}{q'}}\\
&=\sup_{\substack{\|g\|_{\mathcal{S}(I\times\mathbb{R}^{d})}\leq 1}}\left(\int_{I}\left(\int_{\mathbb{R}^{d}}|g(t,x)|^{r}w^{r-2}(x)\ \gamma_{d}(x)\mathrm{d}x\right)^{\frac{q}{r}}\ \mathrm{d}t\right)^{\frac{1}{q}}\\  
&\times\left(\int_{I}\left(\int_{\mathbb{R}^{d}}|f(t,x)|^{r'}w^{r'-2}(x)\gamma_{d}(x)\ \mathrm{d}x\right)^{\frac{q'}{r'}}\ \mathrm{d}t\right)^{\frac{1}{q'}}\\
&\leq\sup_{\substack{\|g\|_{\mathcal{S}(I\times\mathbb{R}^{d})}\leq 1}}\|g\|_{L^{q}_{t}L^{r}_{\gamma_{d}}(w^{r-2})}\cdot \|f\|_{L^{q'}_{t}L^{r'}_{\gamma_{d}}(w^{r'-2})}\\
&\leq \|f\|_{L^{q'}_{t}L^{r'}_{\gamma_{d}}(w^{r'-2})}.
\end{align*}

Therefore, for any $(q,r)\in\mathcal{A}$, we have 
\begin{align}
\label{second}
\|f\|_{\mathcal{N}(I\times \mathbb{R}^{d})}\leq \|f\|_{L^{q'}_{t}L^{r'}_{\gamma_{d}}(w^{r'-2})}.
\end{align}
Before proceeding to the main results of this section, we recall a useful inequality that will be instrumental in estimating the nonlinear term in \eqref{NLS}.
For a more general formulation and detailed proof of this inequality, we refer the reader to \cite[Lemma 3.9]{L. Chaichenets}.

\begin{lemma}
For any \( a, b \in \mathbb{C} \) and \( p \geq 1 \), there exists a constant \( C(p) > 0 \) such that
\begin{equation}
\label{complexinequality}
\left| a|a|^{p-1} - b|b|^{p-1} \right| \leq C(p) |a - b| \big( |a|^{p-1} + |b|^{p-1} \big).
\end{equation}
\end{lemma}

We are now ready to state and prove the main results.

\begin{theorem}[$L^2_{\gamma_d}(\mathbb{R}^d)$ subcritical case]\label{Subcriticalcase}
Let $1 < p < 1 + \frac{4}{d}$ be a subcritical exponent, and let $\mu = \pm 1$. Then the nonlinear Schrödinger equation \eqref{NLS} is locally well-posed in $L^2_{\gamma_d}(\mathbb{R}^d)$ in the subcritical sense.

More precisely, for every $R > 0$, there exists a time $T = T(d, p, R) > 0$ such that for any initial data $u_0 \in B_R := \left\{ u_0 \in L^2_{\gamma_d}(\mathbb{R}^d) : \|u_0\|_{L^2_{\gamma_d}} < R \right\}$, there exists a unique strong $L^{2}_{\gamma_{d}}(\R^{d})$ solution 
\[
u \in \mathcal{S}([-T, T] \times \mathbb{R}^d) \subset C^0_t L^2_{\gamma_d}([-T, T] \times \mathbb{R}^d)
\]
to  \eqref{NLS}.

Furthermore, the solution map $u_0 \mapsto u$ is Lipschitz continuous from $B_R$ to $\mathcal{S}([-T, T] \times \mathbb{R}^d)$.
\end{theorem}

\begin{proof}
Let $I=[-T,T]$, where $0<T\leq\frac{\pi}{2}$ is determined later. Define the Duhamel operator \( D: \mathcal{N}(I \times \mathbb{R}^d) \to \mathcal{S}(I \times \mathbb{R}^d) \) by
\[
(Df)(t,x) := -i \int_0^t e^{-i(t-s)\mathbf{L}} f(s,x) \, \mathrm{d}s.
\]

By employing the Strichartz estimates together with a duality argument, it follows that the operator \( D \) is bounded; more precisely,
\begin{align*}
\|Df\|_{\mathcal{S}(I\times\mathbb{R}^{d})}&=\sup\limits_{(q,r)\in\mathcal{A}}\bigg\|\int_{0}^{t}e^{-i(t-s)\mathbf{L}}f(x,s)\ \mathrm{d}s\bigg\|_{L^{q}_{t}L^{r}_{\gamma_{d}}(w^{r-2})}
\leq C(d)\|f\|_{\mathcal{N}(I\times\mathbb{R}^{d})},
\end{align*}
for some constant $C(d)$ depending on $d$.

 Next, define 
 \[
u_0^*(x,t) := e^{-i t \mathbf{L}} u_0(x).
\]
 Again, by using Strichartz estimate \eqref{HomogeneousStrichartz}
\begin{align*}
\|u_{0}^{*}\|_{\mathcal{S}(I\times\mathbb{R}^{d})}&=\sup\limits_{(q,r)\in\mathcal{A}}\|e^{-it\mathbf{L}}u_{0}\|_{L^{q}_{t}L^{r}_{\gamma_{d}}(w^{r-2})}\\
&\leq \sup\limits_{(q,r)\in\mathcal{A}}C(q,r,d)\|u_{0}\|_{L^{2}_{\gamma_{d}}}.
\end{align*}
Using a uniform constant depending only on \( d \), we conclude
\begin{align}
\|u_{0}^{*}\|_{\mathcal{S}(I\times\mathbb{R}^{d})}\leq C(d) \|u_{0}\|_{L^{2}_{\gamma_{d}}}.
\end{align}

For the iteration scheme, we choose \(\epsilon = 2 C(d) R\) so that the initial function \( u_0^* \) lies in the ball \( B_{\frac{\epsilon}{2}} \).

Next, define the nonlinear operator $N:\mathcal{S}(I\times\mathbb{R}^{d})\to\mathcal{N}(I\times \mathbb{R}^{d})$ as $N(u)=\mu w^{p}|u|^{p-1}u$. We aim to show that for any \( u, v \in B_{\epsilon} \),
\begin{align*}
\|N(u)-N(v)\|_{\mathcal{N}(I\times \mathbb{R}^{d})}\leq \frac{1}{2C(d)}\|u-v\|_{\mathcal{S}(I\times\mathbb{R}^{d})}.
\end{align*}

Using estimate \eqref{second}, for any admissible pair \((q,r) \in \mathcal{A}\), we have
\[
\| N(u) - N(v) \|_{\mathcal{N}(I \times \mathbb{R}^d)} \leq \| N(u) - N(v) \|_{L^{q'}_t L^{r'}_{\gamma_d}(w^{r'-2})} 
= \left\| \| N(u) - N(v) \|_{L^{r'}_{\gamma_d}(w^{r'-2})}\right\|_{L^{q'}_t}.
\]

We begin by estimating \(\|N(u) - N(v)\|_{L^{r'}_{\gamma_d}(w^{r' - 2})}\). By definition, the operator \(N\) satisfies

\[
N(u) = \mu w^{p} |u|^{p-1} u = \mu |u w|^{p-1} (u w).
\]
Setting \(u_1 := u w\), we rewrite \(N(u) = \mu |u_1|^{p-1} u_1\). Similarly, \(N(v) = \mu |v_1|^{p-1} v_1\) where \(v_1 := v w\).

{Applying inequality \eqref{complexinequality}, we obtain
\begin{align}
\nonumber
\|N(u) - N(v)\|_{L^{r'}_{\gamma_d}(w^{r' - 2})} 
&= \| |u_1|^{p-1} u_1 - |v_1|^{p-1} v_1 \|_{L^{r'}_{\gamma_d}(w^{r' - 2})} \\
\label{Firststep}
&\leq C(p) \left\| |u_1 - v_1| \big( |u_1|^{p-1} + |v_1|^{p-1} \big) \right\|_{L^{r'}_{\gamma_d}(w^{r' - 2})},
\end{align}
for some constant \(C(p) > 0\) depending only on \(p\).}

{Consider the pair
\[
(q,r) = \left( \frac{4(p+1)}{d(p-1)},\, p+1 \right).
\]
Since \(1 < p < 1 + \frac{4}{d}\), it follows that
\[
0 < p - 1 < \frac{4}{d}.
\]
Consequently,
\[
\frac{d}{4} < \frac{1}{p-1}.
\]
Therefore,
\[
q = \frac{4(p+1)}{d(p-1)} > \frac{4}{d}(p+1) \cdot \frac{d}{4} = p+1 > 2,
\]
which implies that
\[
2 < q < \infty.
\]
Moreover, we have
\[
r = p + 1 > 2.
\]

Next, we verify the admissibility condition:
\[
\frac{2}{q} + \frac{d}{r} = \frac{d(p-1)}{2(p+1)} + \frac{d}{p+1} = \frac{d(p-1) + 2d}{2(p+1)} = \frac{d(p+1)}{2(p+1)} = \frac{d}{2}.
\]
Thus, the pair \((q,r)\) satisfies the admissibility condition for Strichartz estimates.

Additionally, note that the Hölder conjugate exponent \(r'\) satisfies
\[
r' = \frac{r}{r-1} = \frac{p+1}{p} = \frac{r}{p}.
\]}

For the chosen pair \((q,r)\), from \eqref{Firststep} we have 
\begin{align*}
\|N(u)-N(v)\|_{L^{r'}_{\gamma_{d}}(w^{r'-2})}&\leq C(p)\| |u_{1}-v_{1}|\left(|u_{1}|^{p-1}+|v_{1}|^{p-1}\right)\|_{L^{r'}_{\gamma_{d}}(w^{r'-2})}\\
&=C(p)\pi^{-\frac{d}{2r'}}\| |u_{1}-v_{1}|\left(|u_{1}|^{p-1}+|v_{1}|^{p-1}\right)w\|_{L^{r'}}\\
&\leq C(p)\pi^{\frac{d}{2}\left(\frac{1}{r}-\frac{1}{r'}\right)}\|u_{1}-v_{1}\|_{L^{r'p}_{\gamma_{d}}(w^{r-2})}\| |u_{1}|^{p-1}+|v_{1}|^{p-1}\|_{L^{\frac{pr'}{p-1}}}\\
&\leq C(p)\pi^{\frac{d}{2}\left(\frac{1}{r}-\frac{1}{r'}\right)}\|u_{1}-v_{1}\|_{L^{r}_{\gamma_{d}}(w^{r-2})}\left(\| |u_{1}|^{p-1}\|_{L^{\frac{r}{p-1}}}+ \| |v_{1}|^{p-1}\|_{L^{\frac{r}{p-1}}}\right)\\
&= C(p)\pi^{\frac{d}{2}\left(\frac{1}{r}-\frac{1}{r'}\right)}\|u_{1}-v_{1}\|_{L^{r}_{\gamma_{d}}(w^{r-2})}\left( \|u_{1}\|^{p-1}_{L^{r}}+\|v_{1}\|^{p-1}_{L^{r}}\right)\\
&=C(p)\pi^{\frac{d}{2}\left(\frac{1}{r}-\frac{1}{r'}\right)}\|(u-v)w\|_{L^{r}_{\gamma_{d}}(w^{r-2})}\left( \|uw\|^{p-1}_{L^{r}}+\|vw\|^{p-1}_{L^{r}}\right)\\
&\leq C(p)\pi^{\frac{d}{2}\left(\frac{1}{r}-\frac{1}{r'}+\frac{p-1}{r}\right)}\|u-v\|_{L^{r}_{\gamma_{d}}(w^{r-2})}\left(\|u\|^{p-1}_{L^{r}_{\gamma_{d}}(w^{r-2})}+\|v\|^{p-1}_{L^{r}_{\gamma_{d}}(w^{r-2})}\right)\\
&= C(p)\|u-v\|_{L^{r}_{\gamma_{d}}(w^{r-2})}\left(\|u\|^{p-1}_{L^{r}_{\gamma_{d}}(w^{r-2})}+\|v\|^{p-1}_{L^{r}_{\gamma_{d}}(w^{r-2})}\right),
\end{align*}
where the third inequality follows from H\"older's inequality, the fourth from the triangle inequality and the seventh from the fact that $w(x)\leq 1 \ \forall x\in\mathbb{R}^{d}$.

Using the above estimate and triangle inequality, we obtain
\begin{align}
\nonumber & \,\quad\| N(u)-N(v)\|_{\mathcal{N}(I\times \mathbb{R}^{d})}\\
&\leq \|N(u)-N(v)\|_{L^{q'}_{t}L^{r'}_{\gamma_{d}}(w^{r'-2})}\\
\nonumber &= C(p) \left\| \|N(u)-N(v)\|_{L^{r'}_{\gamma_{d}}(w^{r'-2})}\right\|_{L^{q'}_{t}(I)}\\
\nonumber &\leq C(p)\left\|\|u-v\|_{L^{r}_{\gamma_{d}}(w^{r-2})}\left(\|u\|_{L^{r}_{\gamma_{d}}(w^{r-2})}^{p-1}+\|v\|_{L^{r}_{\gamma_{d}}(w^{r-2})}^{p-1}\right)\right\|_{L^{q'}_{t}(I)}\\
\label{inequlityinbetween}&\leq C(p)\left(\left\|\|u-v\|_{L^{r}_{\gamma_{d}}(w^{r-2})}\|u\|_{L^{r}_{\gamma_{d}}(w^{r-2})}^{p-1}\right\|_{L^{q'}_{t}(I)}            + \left\|\|u-v\|_{L^{r}_{\gamma_{d}}(w^{r-2})}\|v\|_{L^{r}_{\gamma_{d}}(w^{r-2})}^{p-1}\right\|_{L^{q'}_{t}(I)}     \right).
\end{align}
 
Note that
\[
\frac{1}{q'} = \frac{1}{\frac{q}{p}} + \frac{1}{\beta},
\]
where
\[
\frac{1}{\beta} = \frac{1}{q'} - \frac{p}{q} = \frac{d}{4} \left( \frac{4}{d} + 1 - p \right) > 0.
\]
By H\"older's inequality, we have 
\begin{align*}
\left\|\|u-v\|_{L^{r}_{\gamma_{d}}(w^{r-2})}\|u\|_{L^{r}_{\gamma_{d}}(w^{r-2})}^{p-1}\right\|_{L^{q'}_{t}(I)}&\leq \|1\|_{L^\beta_{t}(I)}\left\|\|u-v\|_{L^{r}_{\gamma_{d}}(w^{r-2})}\|u\|_{L^{r}_{\gamma_{d}}(w^{r-2})}^{p-1}\right\|_{L^\frac{q}{p}_{t}(I)}\\
&\leq |I|^{\frac{1}{q'}-\frac{p}{q}}\left\|\|u-v\|_{L^{r}_{\gamma_{d}}(w^{r-2})}\|u\|_{L^{r}_{\gamma_{d}}(w^{r-2})}^{p-1}\right\|_{L^\frac{q}{p}_{t}(I)}.
\end{align*}

Using  similar arguments, we obtain
\[
\left\| \|u - v\|_{L^{r}_{\gamma_d}(w^{r-2})} \|v\|_{L^{r}_{\gamma_d}(w^{r-2})}^{p-1} \right\|_{L^{q'}_t(I)} 
\leq |I|^{\frac{1}{q'} - \frac{p}{q}} \left\| \|u - v\|_{L^{r}_{\gamma_d}(w^{r-2})} \|v\|_{L^{r}_{\gamma_d}(w^{r-2})}^{p-1} \right\|_{L^{\frac{q}{p}}_t(I)}.
\]
Applying these estimates to \eqref{inequlityinbetween}, we obtain
\begin{align}
\label{longinequality}
\nonumber
\| N(u) - N(v) \|_{\mathcal{N}(I\times \mathbb{R}^{d})} 
&\leq C(p) |I|^{\frac{1}{q'} - \frac{p}{q}} \cdot \Bigg( 
\left\| \|u - v\|_{L^{r}_{\gamma_d}(w^{r-2})} \|u\|_{L^{r}_{\gamma_d}(w^{r-2})}^{p-1} \right\|_{L^{\frac{q}{p}}_t(I)} \\
&\quad + \left\| \|u - v\|_{L^{r}_{\gamma_d}(w^{r-2})} \|v\|_{L^{r}_{\gamma_d}(w^{r-2})}^{p-1} \right\|_{L^{\frac{q}{p}}_t(I)} \Bigg).
\end{align}

Since  
\[
\frac{1}{p} + \frac{p-1}{p} = 1,
\]
it follows from H\"older's inequality that
\begin{align}\label{hol1}
\nonumber \left\|\|u-v\|_{L^{r}_{\gamma_{d}}(w^{r-2})}\|u\|_{L^{r}_{\gamma_{d}}(w^{r-2})}^{p-1}\right\|_{L^\frac{q}{p}_{t}(I)}&\leq
\|u-v\|_{L^{q}_{t}L^{r}_{\gamma_{d}}(w^{r-2})}\cdot \|u\|_{L^{q}_{t}L^{r}_{\gamma_{d}}(w^{r-2})}^{p-1}\\
&\leq \|u-v\|_{\mathcal{S}(I\times\mathbb{R}^{d})}\cdot\|u\|_{\mathcal{S}(I\times\mathbb{R}^{d})}^{p-1}.
\end{align}
Similarly,
\begin{align}\label{hol2}
\left\|\|u-v\|_{L^{r}_{\gamma_{d}}(w^{r-2})}\|v\|_{L^{r}_{\gamma_{d}}(w^{r-2})}^{p-1}\right\|_{L^\frac{q}{p}_{t}(I)}
\leq \|u-v\|_{\mathcal{S}(I\times\mathbb{R}^{d})}\cdot \|v\|^{p-1}_{\mathcal{S}(I\times\mathbb{R}^{d})}. 
\end{align}
Using (\ref{hol1}) , (\ref{hol2}) in \eqref{longinequality}, we deduce
\begin{align*}
\| N(u)-N(v)\|_{\mathcal{N}(I\times \mathbb{R}^{d})}\leq C(p)  |I|^{\frac{1}{q'}-\frac{p}{q}}\|u-v\|_{\mathcal{S}(I\times\mathbb{R}^{d})}\left(\|u\|^{p-1}_{\mathcal{S}(I\times\mathbb{R}^{d})}+\|v\|^{p-1}_{\mathcal{S}(I\times\mathbb{R}^{d})}\right).
\end{align*}
Since \(u, v \in B_{\epsilon}\), it follows that
\begin{align}
\label{lastline}
\| N(u)-N(v)\|_{\mathcal{N}(I\times \mathbb{R}^{d})}\leq C(p) 2\epsilon^{p-1}|I|^{\frac{1}{q'}-\frac{p}{q}}\|u-v\|_{\mathcal{S}(I\times\mathbb{R}^{d})}.
\end{align}

Note that the exponent \(\frac{1}{q'} - \frac{p}{q}\) of the factor \(|I|\) is positive. Therefore, by choosing \(T\) sufficiently small, we obtain
\[
\| N(u) - N(v) \|_{\mathcal{N}(I \times \mathbb{R}^d)} \leq \frac{1}{2 C(d)} \| u - v \|_{\mathcal{S}(I \times \mathbb{R}^d)}.
\]

Hence, the iteration principle stated in Proposition \ref{Iterationprinciple} applies, which concludes the proof.
\end{proof}

A crucial aspect in the proof of Theorem \ref{Subcriticalcase} is that, within the subcritical range \(1 < p < 1 + \frac{4}{d}\), the Strichartz estimates provide enhanced control over the nonlinearity. In contrast, for the supercritical case \(p > 1 + \frac{4}{d}\), establishing local existence becomes more involved and may fail. Moreover, the approach employed in Theorem \ref{Subcriticalcase} does not extend to the critical case \(p = 1 + \frac{4}{d}\), since the exponent in \eqref{lastline}, 
\[
\frac{1}{q'} - \frac{p}{q} = \frac{d}{4} \left( \frac{4}{d} + 1 - p \right),
\]
vanishes, resulting in the loss of control over the associated constant.

However, it is still possible to establish existence in the critical sense. We have the following result.

{\begin{theorem}[Critical case]\label{criticalcase}
Let $p = 1 + \frac{4}{d}$ be the $L^2_{\gamma_d}$-critical exponent, and let $\mu = \pm 1$. Then the nonlinear Schrödinger equation \eqref{NLS} is locally well-posed in $L^2_{\gamma_d}(\mathbb{R}^d)$ in the critical sense.

More precisely, given any $R > 0$, there exists $\eta = \eta(R, d) > 0$ such that for every $u_1 \in L^2_{\gamma_d}(\mathbb{R}^d)$ with $\|u_1\|_{L^2_{\gamma_d}} \leq R$, and for any time interval $I = I(\eta) \subseteq \left[-\frac{\pi}{2}, \frac{\pi}{2}\right]$ containing $0$ on which the linear evolution satisfies the smallness condition
\begin{align}
\label{existenceinterval}
\|e^{-it\mathbf{L}}u_1\|_{L^{p+1}_t(I; L^{p+1}_{\gamma_d}(w^{p-1}))} \leq \eta,
\end{align}
then for all initial data $u_0 \in B(u_1, \eta) := \left\{ u_0 \in L^2_{\gamma_d}(\mathbb{R}^d) : \|u_0 - u_1\|_{L^2_{\gamma_d}} \leq \eta \right\}$, there exists a unique strong  $L^{2}_{\gamma_{d}}(\mathbb{R}^{d})$ solution $u \in \mathcal{S}(I \times \mathbb{R}^d)$ to the integral equation \eqref{Duhamel}.

Furthermore, the solution map $u_0 \mapsto u$ is Lipschitz continuous from $B(u_1, \eta)$ to $\mathcal{S}(I \times \mathbb{R}^d)$, with Lipschitz constant depending on $u_1$.
\end{theorem}}

Before proving the above theorem, we verify that the condition \eqref{existenceinterval} is meaningful.

\begin{lemma}
Let \(p = 1 + \frac{4}{d}\) and \(u_1 \in L^2_{\gamma_d}\) satisfies \(\|u_1\|_{L^2_{\gamma_d}} \leq R\). Then for any \(\eta > 0\), there exists a time interval \(I = I(\eta) \subseteq \left[-\frac{\pi}{2}, \frac{\pi}{2}\right]\) containing \(0\) such that
\[
\| e^{-i t \mathbf{L}} u_1 \|_{L^{p+1}_t(I; L^{p+1}_{\gamma_d}(w^{p-1}))} \leq \eta.
\]
\end{lemma}

{\begin{proof}
Consider the sequence \(\{g_n(t)\}_{n \in \mathbb{N}}\), where
\[
g_n(t) = \| e^{-i t \mathbf{L}} u_1 \|_{L^{p+1}_{\gamma_d}(w^{p-1})} \, \chi_{I_n}(t), \quad I_n = \left[-\frac{1}{n}, \frac{1}{n}\right], \quad \forall n \in \mathbb{N}.
\]
Note that \(g_n(t) \to 0\) almost everywhere in \(\left[-\frac{\pi}{2}, \frac{\pi}{2}\right]\), and
\[
\| g_n \|_{L^{p+1}_t}^{p+1} = \int_{-\frac{\pi}{2}}^{\frac{\pi}{2}} |g_n(t)|^{p+1} \, dt = \| e^{-i t \mathbf{L}} u_1 \|_{L^{p+1}_t(I_n; L^{p+1}_{\gamma_d}(w^{p-1}))}^{p+1}.
\]
Moreover, the function \(|g_n(t)|^{p+1}\) is dominated by \(\| e^{-i t \mathbf{L}} u_1 \|_{L^{p+1}_{\gamma_d}(w^{p-1})}^{p+1}\)and this function belongs to  \(L^1_t\left[-\frac{\pi}{2}, \frac{\pi}{2}\right]\)  due to the  homogeneous Strichartz estimate \eqref{HomogeneousStrichartz}.
Hence, by the Dominated Convergence Theorem, 
\[
\| e^{-i t \mathbf{L}} u_1 \|_{L^{p+1}_t(I_n; L^{p+1}_{\gamma_d}(w^{p-1}))}^{p+1} \to 0
\]
as \(n \to \infty\).

Therefore, for any \(\eta > 0\), there exists \(N = N(\eta) \in \mathbb{N}\) such that for all \(n \geq N\),
\[
\| e^{-i t \mathbf{L}} u_1 \|_{L^{p+1}_t(I_n; L^{p+1}_{\gamma_d}(w^{p-1}))} \leq \eta.
\]
\end{proof}}

{\begin{proof}[Proof of Theorem \ref{criticalcase}]
As in the subcritical case, we apply the iteration principle (Proposition \ref{Iterationprinciple}). The iteration is carried out in the same function spaces $\mathcal{N}(I \times \mathbb{R}^d)$ and $\mathcal{S}(I \times \mathbb{R}^d)$, and uses the same operators $D$ and $N$ as in the proof of Theorem \ref{Subcriticalcase}. However, we now introduce a modified space $\mathcal{S}^{0}(I \times \mathbb{R}^d)$, equipped with the norm
\[
\|u\|_{\mathcal{S}^{0}(I \times \mathbb{R}^d)} := \delta \|u\|_{\mathcal{S}(I \times \mathbb{R}^d)} + \|u\|_{L^{p+1}_t L^{p+1}_{\gamma_d}(w^{p-1})},
\]
where $0 < \delta < 1$ is a parameter to be chosen later. Let us also denote by $u_0^* := e^{-it \mathbf{L}} u_0$  corresponding to initial data $u_0 \in B(u_1, \eta)$, where $\|u_0 - u_1\|_{L^2_{\gamma_d}} \leq \eta$.

We first estimate the norm of $u_0^*$ in $\mathcal{S}^{0}$:
\begin{align*}
\|u_0^*\|_{\mathcal{S}^{0}(I \times \mathbb{R}^d)} 
&= \delta \|u_0^*\|_{\mathcal{S}(I \times \mathbb{R}^d)} + \|u_0^*\|_{L^{p+1}_t L^{p+1}_{\gamma_d}(w^{p-1})} \\
&\leq \delta C(d) \|u_0\|_{L^2_{\gamma_d}} + \|e^{-it\mathbf{L}}(u_0 - u_1)\|_{L^{p+1}_t L^{p+1}_{\gamma_d}(w^{p-1})} + \|e^{-it\mathbf{L}} u_1\|_{L^{p+1}_t L^{p+1}_{\gamma_d}(w^{p-1})} \\
&\leq \delta C(d)(R + \eta) + \|u_0 - u_1\|_{L^2_{\gamma_d}} + \eta \\
&\leq \delta C(d)(R + \eta) + 2\eta =: \frac{\epsilon}{2},
\end{align*}
where we have used the Strichartz estimate and triangle inequality.

Next, we show that the operator $D$ is bounded from $\mathcal{N}$ to $\mathcal{S}^{0}$. From the homogeneous Strichartz estimate, we have
\[
\|Df\|_{\mathcal{S}(I \times \mathbb{R}^d)} \leq C(d) \|f\|_{\mathcal{N}(I \times \mathbb{R}^d)}.
\]
Moreover, since $(p+1, p+1)$ is an admissible pair and $\mathcal{S}(I \times \mathbb{R}^d) \subseteq L^{p+1}_t L^{p+1}_{\gamma_d}(w^{p-1})$, we get
\[
\|Df\|_{L^{p+1}_t L^{p+1}_{\gamma_d}(w^{p-1})} \leq \|Df\|_{\mathcal{S}(I \times \mathbb{R}^d)} \leq C(d) \|f\|_{\mathcal{N}(I \times \mathbb{R}^d)}.
\]
Therefore,
\begin{align*}
\|Df\|_{\mathcal{S}^{0}(I \times \mathbb{R}^d)} 
&= \delta \|Df\|_{\mathcal{S}} + \|Df\|_{L^{p+1}_t L^{p+1}_{\gamma_d}(w^{p-1})} \\
&\leq \delta C(d) \|f\|_{\mathcal{N}} + C(d) \|f\|_{\mathcal{N}} \\
&\leq C(d)(1 + \delta) \|f\|_{\mathcal{N}} \leq 2C(d) \|f\|_{\mathcal{N}},
\end{align*}
where we used that $0 < \delta < 1$.

Now, using the same arguments as in the subcritical case, we estimate the nonlinearity:
\begin{align*}
\|N(u) - N(v)\|_{\mathcal{N}(I \times \mathbb{R}^d)} 
&\leq C(d) |I|^{\frac{1}{q'} - \frac{p}{q}} \|u - v\|_{L^{p+1}_t L^{p+1}_{\gamma_d}(w^{p-1})} \\
&\quad \times \left( \|u\|^{p+1}_{L^{p+1}_t L^{p+1}_{\gamma_d}(w^{p-1})} + \|v\|^{p+1}_{L^{p+1}_t L^{p+1}_{\gamma_d}(w^{p-1})} \right).
\end{align*}
Since we are in the critical case \(p = 1 + \frac{4}{d}\), we observe that
\[
\frac{1}{q'} - \frac{p}{q} = \frac{d}{4} \left( \frac{4}{d} + 1 - p \right) = 0.
\]
Thus, the time factor \( |I|^{\frac{1}{q'} - \frac{p}{q}} \) does not provide smallness. However, since \(u, v \in B_\epsilon\), we obtain
\[
\|N(u) - N(v)\|_{\mathcal{N}(I \times \mathbb{R}^d)} \leq 2 C(d) \epsilon^{p+1} \|u - v\|_{\mathcal{S}^{0}(I\times \mathbb{R}^{d})}.
\]
Now, since \(\epsilon\) depends on \(\eta\) and \(\delta\), we can choose them sufficiently small so that
\[
\|N(u) - N(v)\|_{\mathcal{N}(I \times \mathbb{R}^d)} \leq \frac{1}{4C(d)} \|u - v\|_{\mathcal{S}^{0}(I\times \mathbb{R}^{d})}.
\]
Combining all of the above, we can apply the iteration principle to obtain a unique solution \(u \in \mathcal{S}^{0}(I \times \mathbb{R}^d)\), which in turn belongs to \(\mathcal{S}(I \times \mathbb{R}^d)\), completing the proof.
\end{proof}
}


\end{document}